\documentclass[12pt]{article}
\usepackage{amssymb,amsfonts,amsmath,amsthm,cite}
\usepackage{epsfig}
\newtheorem{thm}{Theorem}[section]
\newtheorem{cor}{Corollary}[section]
\newtheorem{lem}{Lemma}[section]

\makeatletter \@addtoreset{equation}{section}

\def\pf{\noindent {\it Proof.\ }}
\def\qed{\hfill \rule{4pt}{7pt}}

\begin{document}

\begin{center}
{\large {\bf Arithmetic Properties of Overpartition Pairs }}

\vskip 6mm

{\small William Y.C. Chen$^1$ and Bernard L.S. Lin$^2$
\\[%
2mm] Center for Combinatorics, LPMC-TJKLC\\
Nankai University, Tianjin 300071,
P.R. China \\[3mm]
$^1$chen@nankai.edu.cn, $^2$linlishuang@cfc.nankai.edu.cn \\[0pt%
] }
\end{center}

\noindent {\bf Abstract.} Bringmann and Lovejoy
 introduced a rank for overpartition pairs and
investigated its role in congruence properties of
$\overline{pp}(n)$, the number of overpartition pairs of $n$.
 In particular, they applied
 the theory of Klein forms to show that there
 exist many Ramanujan-type congruences for
the number $\overline{pp}(n)$. In this paper, we shall derive two
Ramanujan-type identities and some explicit congruences for
$\overline{pp}(n)$. Moreover, we
find three ranks as combinatorial interpretations
of the fact that $\overline{pp}(n)$ is divisible by three
for any $n$.  We also
 construct  infinite families of congruences for
$\overline{pp}(n)$ modulo $3$, $5$,
and $9$.





\section{Introduction}

A partition of a positive integer $n$ is
a non-increasing sequence
of positive integers whose sum is $n$.
An overpartition $\lambda$ of
$n$ is a partition of $n$ for which
 the first occurrence of a number
may be overlined. Let $\overline{p}(n)$ denote the number of
overpartitions of $n$. Congruence properties for $\overline{p}(n)$
have been extensively studied, see, for example, Fortin, Jacob and
Mathieu \cite{FJM05}, Hirschhorn and Sellers \cite{HS05b}, Kim
\cite{Kim09}, Lovejoy and Osburn \cite{Lovejoy10}, and Mahlburg
\cite{Mahlburg04}. In this paper, we shall be concerned with
arithmetic properties of the number of overpartition pairs of $n$.
 An
overpartition pair $\pi$ of $n$ is a pair of overpartitions
$(\lambda,\mu)$ such that the sum of all of the parts is $n$. Note
that either $\lambda$ or $\mu$ may be an overpartition of zero.  We
need to pay special attention to the overpartition of zero. There is
only one partition of zero, and there is only one overpartition of
zero as well. Let $\overline{pp}(n)$ denote the number of
overpartition pairs of $n$.
 Then the generating function for
$\overline{pp}(n)$ is
\begin{equation}\label{eqgenpp1}
\sum_{n=0}^\infty \overline{pp}(n)q^n= \frac{(-q;q)_\infty
^2}{(q;q)_\infty ^2}.
\end{equation}
Throughout this paper, we adopt the following standard
$q$-series notation for $|q|<1$,
\[
(a;q)_\infty =\prod_{k=1}^\infty (1-aq^{k-1}).
\]

 Bringmann and Lovejoy \cite{Lovejoy08} defined
 a rank for overpartition pairs
to investigate
congruence properties of $\overline{pp}(n)$.
 Let
$\overline{NN}(m,n)$ denote the number of overpartition pairs of $n$
with rank $m$, and let $\overline{NN}(r,t,n)$ denote
 the number of overpartition pairs
of $n$ with rank congruent to $r$ modulo $t$. They obtained a bivariate
 generating function for
$\overline{NN}(m,n)$ from which they
derived the following relation  for $0\leq r \leq 2$,
\[
\overline{NN}(r,3,3n+2)=\frac{\overline{pp}(3n+2)}{3}.
\]
This leads to the following Ramanujan-type congruence
\begin{equation}\label{eqcong3}
\overline{pp}(3n+2) \equiv 0\ ({\rm mod}\ 3).
\end{equation}
Furthermore, by using the theory of Klein forms,
Bringmann and Lovejoy
 \cite{Lovejoy08} proved that there exist
  infinitely many  Ramanujan-type congruences
  for $\overline{pp}(n)$.
  Let $l$ be an
odd prime and let $t$ be an odd number which is a power
of $l$ or relatively prime to $l$.
Then for any positive integer $j$, there
are infinitely many non-nested arithmetic progressions
$An + B$ such that
\begin{equation}\label{eqrankcong}
\overline{NN}(r,t,An+B)\equiv 0\ ({\rm mod}\ l^j)
\end{equation}
for any $0\leq r\leq t-1$.
Hence there are infinitely many non-nested
arithmetic progressions $An +
B$ satisfying
\begin{equation}\label{eqcong}
 \overline{pp}(An+B)\equiv 0\ ({\rm mod}\ l^j)
\end{equation}
for any odd prime $l$ and any positive integer $j$. For the case
$l=2$, using the theory of modular forms, they have shown that
\eqref{eqcong} holds for any positive integer $j$.

However,  the theory of Klein forms used to derive the congruence
relation \eqref{eqcong} is not constructive and it does not give
explicit  arithmetic progressions $An+B$ in the statement. So it
is still desirable to
 find
explicit congruences for $\overline{pp}(n)$.
In this paper, we obtain some congruences
 modulo $3$ and $5$.

For the case of modulo $3$, we obtain a Ramanujan-type
identity
\begin{equation}
\sum_{n=0}^\infty \overline{pp}(3n+2)q^n=12\frac{(q^2;q^2)_\infty
^6(q^3;q^3)_\infty ^6}{(q;q)_\infty ^{14}},
\end{equation} which implies
\eqref{eqcong3}.
Furthermore, we show that there are  infinite families of
congruences modulo $3$ satisfied by $\overline{pp}(n)$.
For example,  for
any $\alpha \geq 1$ and $n\geq 0$,
\begin{equation}
\overline{pp}(9^\alpha(3n+1))\equiv\overline{pp}(9^\alpha(3n+2))\equiv
0\ ({\rm mod}\ 3).
\end{equation}

For the case of modulo $5$, we obtain three
Ramanujan-type congruences
\begin{equation}
\overline{pp}(20n+11)\equiv\overline{pp}(20n+15)
\equiv\overline{pp}(20n+19)\equiv 0\ ({\rm mod}\ 5),
\end{equation} for any  $n\geq 0$.
We also obtain infinite families of congruences modulo $5$.
For example, for any $\alpha \geq 1$ and $n\geq 0$,
\begin{equation}
\overline{pp}(5^\alpha(5n+2))\equiv
\overline{pp}(5^\alpha(5n+3))\equiv 0\ ({\rm mod}\ 5).
\end{equation}

Motivated by the
work of Paule and Radu \cite{Paule09} on some strange
congruences in their words, we
establish similar  congruences for $\overline{pp}(n)$.
For example,
for any $k\geq 0$,
\begin{equation}\label{eqsample2}
 \overline{pp}(5\cdot29^k)\equiv 3(k+1)\ ({\rm mod\ }5)
\end{equation}
and
\begin{equation}\label{eqsample1}
 \overline{pp}(2\cdot13^k)\equiv 3(k+1)\ ({\rm mod\ }9).
\end{equation}

In order to give combinatorial interpretations of the fact that
$\overline{pp}(3n+2)$ is divisible by $3$ for any $n\geq0$, we find
three ranks of overpartition pairs that serve this purpose.

This paper is organized as follows. In Section 2, we obtain two
Ramanujan-type identities  and some Ramanujan-type congruences
modulo $5$ and $64$. In Section 3, we give three combinatorial
interpretations for the congruence \eqref{eqcong3}. Section 4
contains infinite families of congruences modulo $3$ and $5$. In
Section 5 is concerned with congruences modulo $5$ and $9$ in the
flavor of the strange congruences of Paule and Radu.

\section{Ramanujan-type identities and congruences }

In this section,  we  establish
 two Ramanujan-type
identities and derive some congruence relations modulo $5$ and $64$.

\begin{thm}\label{thmRamtypeiden} We have
\begin{eqnarray}
\sum_{n=0}^\infty \overline{pp}(3n+2)q^n&=&12\frac{(q^2;q^2)_\infty
^6(q^3;q^3)_\infty ^6}{(q;q)_\infty ^{14}},\label{eq3n+2}\\[5pt]
\sum_{n=0}^\infty \overline{pp}(4n+3)q^n&=&32\frac{(q^2;q^2)_\infty
^{20}}{(q;q)_\infty ^{22}}.\label{eq4n+3}
\end{eqnarray}
\end{thm}

To prove the above identities,
 we recall  two
Ramanujan's theta functions $\varphi(q)$ and $\psi(q)$, namely,
\[
\varphi(q)=\sum_{n=-\infty}^\infty q^{n^2},\quad
\psi(q)=\sum_{n=0}^\infty q^{n(n+1)/2}.
\]
The following two identities are due to Gauss, see, for example,
Berndt \cite[p.11]{Berndt06}.
\[
\varphi(-q)=\frac{(q;q)_\infty^2}{(q^2;q^2)_\infty},\quad
\psi(q)=\frac{(q^2;q^2)_\infty^2}{(q;q)_\infty}.
\]
As shown by Hirschhorn and Sellers \cite{HS05}, the the generating
function of $\overline{p}(n)$ is
\[
\sum_{n=0}^\infty \overline{p}(n)q^n=\frac{1}{\varphi(-q)}.
\]
This implies that the generating function
of $\overline{pp}(n)$ equals
\begin{equation}\label{eqgenpp2}
\sum_{n=0}^\infty \overline{pp}(n)q^n=\frac{1}{\varphi(-q)^2}.
\end{equation}

The following dissection formula of Hirschhorn and Sellers
\cite{HS05} plays a key role
 in  the proof of Theorem
\ref{thmRamtypeiden}.

\begin{lem}\label{lembyHS}
\begin{eqnarray}
\frac{1}{\varphi(-q)}
&=&\frac{\varphi(-q^9)}{\varphi(-q^3)^4}\left(\varphi(-q^9)^2+2q\varphi(-q^9)A(q^3)+4q^2
A(q^3)^2\right)\label{eqlem3n}\\[5pt]
&=&\frac{1}{\varphi(-q^4)^4}
\left(\varphi(q^4)^3+2q\varphi(q^4)^2\psi(q^8)+4q^2\varphi(q^4)
\psi(q^8)^2+8q^3\psi(q^8)^3\right),\label{eqlem4n}
\end{eqnarray}
where
\[
A(q)=\frac{(q;q)_\infty(q^6;q^6)_\infty^2}
{(q^2;q^2)_\infty(q^3;q^3)_\infty}.
\]
\end{lem}

 \noindent {\it Proof of Theorem \ref{thmRamtypeiden}.}
  Applying the $3$-dissection formula \eqref{eqlem3n} in \eqref{eqgenpp2}, we see that
\begin{equation}\label{eqHS1}
\sum_{n=0}^\infty \overline{pp}(n)q^n =
\frac{\varphi(-q^9)^2}{\varphi(-q^3)^8}\left(\varphi(-q^9)^2+2q\varphi(-q^9)A(q^3)+4q^2
A(q^3)^2\right)^2.
\end{equation}
Choosing  those terms on each side of the above identity for which
the powers of $q$ are of the form $3n+2$ and noting that all the
series \eqref{eqHS1} are functions of  $q^3$ if the factors $q$ and
$q^2$ are not taken into account, we find that
\begin{eqnarray*}
\sum_{n=0}^\infty \overline{pp}(3n+2)q^{3n+2}&=&
\frac{\varphi(-q^9)^2}{\varphi(-q^3)^8}\left(8q^2\varphi(-q^9)^2A(q^3)^2+4q^2\varphi(-q^9)^2A(q^3)^2\right)\\[5pt]
&=&12q^2A(q^3)^2\frac{\varphi(-q^9)^4}{\varphi(-q^3)^8}.
\end{eqnarray*}
 Dividing both sides of the above identity by
$q^2$ and replacing $q^3$ by $q$, we obtain that
\[
\sum_{n=0}^\infty \overline{pp}(3n+2)q^n
=12A(q)^2\frac{\varphi(-q^3)^4}{\varphi(-q)^8}.
\]
This yields \eqref{eq3n+2}. Similarly,
\[
\sum_{n=0}^\infty \overline{pp}(n)q^n =
\frac{1}{\varphi(-q^4)^8}\left(\varphi(q^4)^3+2q\varphi(q^4)^2\psi(q^8)+4q^2\varphi(q^4)
\psi(q^8)^2+8q^3\psi(q^8)^3\right)^2.
\]
Choosing the terms in the above identity for which the
powers of $q$ are of the form $4n+3$, we find that
\begin{eqnarray*}
\sum_{n=0}^\infty
\overline{pp}(4n+3)q^{4n+3}&=&\frac{1}{\varphi(-q^4)^8}\left(16q^3\varphi(q^4)^3\psi(q^8)^3+16q^3\varphi(q^4)^3\psi(q^8)^3\right)\\[5pt]
&=&32q^3\frac{\varphi(q^4)^3\psi(q^8)^3}{\varphi(-q^4)^8}.
\end{eqnarray*}
Dividing both sides of the above identity  by $q^3$
and replacing
$q^4$ by $q$, we deduce that
\begin{equation}\label{eqpsi4n+3}
\sum_{n=0}^\infty
\overline{pp}(4n+3)q^n=32\frac{\varphi(q)^3\psi(q^2)^3}{\varphi(-q)^8},
\end{equation}
which is equivalent to \eqref{eq4n+3}.
 This completes the proof. \qed

In view of Theorem \ref{thmRamtypeiden}, it can be seen that
$\overline{pp}(3n+2)$ and $\overline{pp}(4n+3)$ are divisible by
$4$. In fact, for all $n\geq 1$, $\overline{pp}(n)$ is divisible
by $4$, since
\begin{eqnarray*}
\sum_{n=0}^\infty \overline{pp}(n)q^n&\equiv&
\left(1+2\sum_{n=0}^\infty (-q)^{n^2}\right)^2\sum_{n=0}^\infty
\overline{pp}(n)q^n \ ({\rm mod\
}4)\\[5pt]
&=&\varphi(-q)^2\frac{1}{\varphi(-q)^2}=1.
\end{eqnarray*}
Indeed,  Keister, Sellers and Vary \cite{KSV09} have shown
 that for $n\geq 1$,
\[
\overline{pp}(n)\equiv \begin{cases} 4\ ({\rm mod\ }8), & \text{ if
}n \text{ is a square or twice a square},\\ 0\ ({\rm mod\ }8), &
\text{ otherwise.}
\end{cases}
\]

With the aid of  \eqref{eq4n+3} and
 the following  relations for any prime $p$,
\begin{equation}\label{eqprod1}
(q;q)_\infty^p \equiv (q^p;q^p)_\infty\ ({\rm mod\ }p),
\end{equation}
we are led to the following congruence relations modulo
  $5$ and $64$.

\begin{cor}\label{thmRamtypecong} For any
nonnegative integer $n$,
\begin{eqnarray}
\overline{pp}(8n+7)&\equiv& 0\ ({\rm mod}\ 64),\label{eqcong8n+7}\\[5pt]
\overline{pp}(20n+11)&\equiv& 0\ ({\rm mod}\ 5),\label{eqcong20n+11}\\[5pt]
\overline{pp}(20n+15)&\equiv& 0\ ({\rm mod}\ 5),\label{eqcong20n+15}\\[5pt]
\overline{pp}(20n+19)&\equiv& 0\ ({\rm mod}\
5).\label{eqcong20n+19}
\end{eqnarray}
\end{cor}

\pf From \eqref{eq4n+3} and  \eqref{eqprod1} with $p=2$, we have
\[
\sum_{n=0}^\infty \frac{\overline{pp}(4n+3)}{32}q^n\equiv
\frac{(q^2;q^2)_\infty ^{20}}{(q^2;q^2)_\infty ^{11}}\equiv
(q^2;q^2)_\infty^9 \ ({\rm mod\ }2).
\]
 This
yields congruence \eqref{eqcong8n+7} by equating the coefficients of
$q^{2n+1}$ for $n\geq 0$. Again by \eqref{eq4n+3} and
\eqref{eqprod1} with $p=5$, we see that
\begin{equation}\label{eqtempthmcong1}
\sum_{n=0}^\infty \overline{pp}(4n+3)q^n\equiv
2\frac{(q^{10};q^{10})_\infty^4}{(q^5;q^5)_\infty ^4}\cdot
\frac{1}{(q;q)_\infty ^2}\ ({\rm mod\ }5).
\end{equation}
Let $p_{-2}(n)$ be defined by
\[
\sum_{n=0}^\infty p_{-2}(n)q^n = \frac{1}{(q;q)_\infty^2}.
\]
It has been shown by Ramanathan \cite{Ramanathan50}
that  for $n\geq 0$,
\[
p_{-2}(5n+2)\equiv p_{-2}(5n+3)\equiv p_{-2}(5n+4)\equiv 0\ ({\rm
mod\ }5).
\]
Combining \eqref{eqtempthmcong1} and the above three
 congruences, we
deduce the  congruence relations \eqref{eqcong20n+11},
\eqref{eqcong20n+15} and \eqref{eqcong20n+19}.
 This completes the proof. \qed

\section{Three ranks for overpartition pairs}

In this section, we give three combinatorial interpretations for the
fact that $\overline{pp}(3n + 2)$ is divisible by $3$.

The first rank of an overpartition pair
 $\pi=(\lambda,\mu\mathbb{})$, denoted $r_1(\pi)$,
is defined to be
$n_1(\lambda)-n_1(\mu)$, where $n_1(\lambda)$ denotes
 the number of
parts of an overpartition $\lambda$. As usual,
let $R_1(m,n)$ denote the number of
overpartition pairs of $n$ with   $r_1(\pi)=m$ and let
$R_1(s,t,n)$ denote the number of overpartition pairs of $n$
with $r_1(\pi)\equiv s\ ({\rm mod\ }t)$. By symmetry, we see that
$R_1(m,n)=R_1(-m,n)$, and so $R_1(s,t,n)=R_1(t-s,t,n)$.
It is easy to
derive the bivariate generating function for $R_1(m,n)$,
that is,
\begin{equation}\label{eqrank1}
\sum_{m=-\infty}^\infty \sum_{n=0}^\infty
R_1(m,n)z^mq^n=\frac{(-qz;q)_\infty}{(qz;q)_\infty}\cdot
\frac{(-q/z;q)_\infty}{(q/z;q)_\infty}.
\end{equation}
Here we adopt the convention that the empty overpartition pair of
$0$ has rank zero and this convention also holds for the other
 two ranks that will be introduced in this section.
The following theorem  shows that the rank $r_1(\pi)$
 leads to a classification of
overpartition pairs of $3n+2$ into three equinumerous sets.

\begin{thm}For  $0\leq s\leq 2$, we have
\begin{equation}\label{eqrank1a}
R_1(s,3,3n+2)=\frac{\overline{pp}(3n+2)}{3}.
\end{equation}
\end{thm}
\pf  Substituting $z=\xi=e^{2\pi i/3}$ into \eqref{eqrank1} and
using the symmetry relation $R_1(1,3,n)=R_1(2,3,n)$,
 we find that
\begin{eqnarray}
\sum_{n=0}^\infty
(R_1(0,3,n)-R_1(1,3,n))q^n&=&\frac{(-q\xi;q)_\infty(-q\xi^2;q)_\infty}{(q\xi;q)_\infty(q\xi^2;q)_\infty}\nonumber\\[5pt]
&=&\frac{(-q^3;q^3)_\infty}{(q^3;q^3)_\infty}\cdot
\frac{(q;q)_\infty}{(-q;q)_\infty}\nonumber\\[5pt]
&=&\frac{(-q^3;q^3)_\infty}{(q^3;q^3)_\infty}
\sum_{n=-\infty}^\infty (-1)^nq^{n^2}.\label{eqrank1b}
\end{eqnarray}
Here the second equality
follows from identity
\[
(1-x^3)=(1-x)(1-x\xi)(1-x\xi^2).
\] Equating the coefficients
of $q^{3n+2}$ on both sides of
 \eqref{eqrank1b}, and
observing that there are no squares congruent to $2$ modulo
$3$, we conclude that
\[
R_1(0,3,3n+2)=R_1(1,3,3n+2),
\]
and so
\[
R_1(0,3,3n+2)=R_1(1,3,3n+2)=R_1(2,3,3n+2)=\frac{\overline{pp}(3n+2)}{3}.
\]
This completes the proof.\qed

We now give the second rank $r_2$. Let
$\pi=(\lambda,\mu)$ be an overpartition pair.
Define
$r_2(\pi)=n_2(\lambda)-n_2(\mu)$, where $n_2(\lambda)$
denotes the
number of overlined parts of an overpartition
$\lambda$. Similarly, let $R_2(m,n)$
denote the number of overpartition pairs of $n$ with
$r_2(\pi)=m$ and let $R_2(s,t,n)$ denote
the number of overpartition
pairs of $n$ with $r_2(\pi)\equiv s\ ({\rm mod\ }t)$.
Then we have the following relation.

\begin{thm} For $n\geq 0$, we have
\begin{equation}\label{eqrank2}
R_2(0,3,3n+2)\equiv R_2(1,3,3n+2)\equiv R_2(2,3,3n+2)\ ({\rm mod\
}3).
\end{equation}
\end{thm}

\pf It is routine to check that
\begin{equation}\label{eqrank2a}
\sum_{m=-\infty}^\infty \sum_{n=0}^\infty
R_2(m,n)z^mq^n=\frac{(-qz;q)_\infty}{(q;q)_\infty} \cdot
\frac{(-q/z;q)_\infty}{(q;q)_\infty}.
\end{equation}
Using the fact that $R_2(1,3,n)=R_2(2,3,n)$ and setting
$z=\xi=e^{2\pi i/3}$ in \eqref{eqrank2a}, we find
\begin{eqnarray} \sum_{n=0}^\infty
(R_2(0,3,n)-R_2(1,3,n))q^n&=&\frac{(-q\xi;q)_\infty(-q\xi^2;q)_\infty}{(q;q)_\infty^2}\nonumber\\[5pt]
&=&\frac{(-q^3;q^3)_\infty}{(q;q)_\infty(q^2;q^2)_\infty}\label{eqrank2b}\\[5pt]
&=&\frac{(-q^3;q^3)_\infty}{(q;q)_\infty^3}\sum_{n=-\infty}^\infty
(-q)^{n^2}\nonumber\\[5pt]
&\equiv&\frac{(-q^3;q^3)_\infty}{(q^3;q^3)_\infty}\sum_{n=-\infty}^\infty
(-q)^{n^2}\ ({\rm mod\ }3).\nonumber
\end{eqnarray}
Since there are no squares congruent to $2$ modulo $3$, we
see that
\[
R_2(0,3,3n+2)-R_2(1,3,3n+2)\equiv 0\ ({\rm mod\ }3),
\]
 and hence
the proof is complete.\qed

It is worth mentioning that
 Andrews, Lewis and Lovejoy \cite{ALL02}
 investigated the
arithmetic properties of the number $PD(n)$
 of partitions of $n$
with designated summands, whose generating function is given by
\eqref{eqrank2b}, that is,
\[
\sum_{n=0}^\infty PD(n)q^n=
\frac{(q^6;q^6)_\infty}{(q;q)_\infty(q^2;q^2)_\infty
(q^3;q^3)_\infty}.
\]
For example, it has been shown that $PD(3n+2)$  is divisible by
three. It should also be mentioned that Chan \cite{Chan08} studied
the number $a(n)$ given by
\[
\sum_{n=0}^\infty a(n)q^n
=\frac{1}{(q;q)_\infty(q^2;q^2)_\infty},
\]
and derived a Ramanujan-type identity for
$a(3n+2)$, that is,
\begin{equation}\label{eqChan}
\sum_{n=0}^\infty
a(3n+2)q^n=3\frac{(q^3;q^3)_\infty^3(q^6;q^6)_\infty^3}{(q;q)_\infty^4(q^2;q^2)_\infty^4}.
\end{equation}
From  \eqref{eqrank2b} and \eqref{eqChan}, we get the following
formula.

\begin{cor} We have
\begin{equation}
\sum_{n=0}^\infty \left(R_2(0,3,3n+2)-R_2(1,3,3n+2)\right)q^n=
3\frac{(q^3;q^3)_\infty^3(q^6;q^6)_\infty^3}{(q;q)_\infty^5(q^2;q^2)_\infty^3}.
\end{equation}
\end{cor}

Finally, we turn to  the third rank $r_3$ of an overpartition pair
$\pi=(\lambda,\mu)$, which is defined by
$r_3(\pi)=n_3(\lambda)-n_3(\mu)$, where $n_3(\lambda)$ denotes the
number of non-overlined parts of an overpartition $\lambda$.
Similarly, let $R_3(m,n)$ denote the number of overpartition pairs
of $n$ with  $r_3(\pi)=m$ and let $R_3(s,t,n)$ denote the number of
overpartition pairs of $n$ with
 $r_3(\pi)\equiv s\ ({\rm
mod\ }t)$. Then we have the following theorem.

\begin{thm}
For $0\leq s\leq 2$, we have
\begin{equation}\label{eqrank3}
R_3(s,3,3n+2)=\frac{\overline{pp}(3n+2)}{3}.
\end{equation}
\end{thm}

\pf It is easy to derive that
\begin{equation}\label{eqrank3a}
\sum_{m=-\infty}^\infty \sum_{n=0}^\infty
R_3(m,n)z^mq^n=\frac{(-q;q)_\infty^2}{(qz;q)_\infty(q/z;q)_\infty}.
\end{equation}
Using the fact that $R_3(1,3,n)=R_3(2,3,n)$ and setting
$z=\xi=e^{2\pi i/3}$ in \eqref{eqrank3a}, we find that
\begin{eqnarray*}
\sum_{n=0}^\infty
\left(R_3(0,3,n)-R_3(1,3,n)\right)q^n&=&\frac{(-q;q)_\infty^2}{(q\xi;q)_\infty(q/\xi;q)_\infty}\\[5pt]
&=&\frac{(-q;q)_\infty^2(q;q)_\infty}{(q^3;q^3)_\infty}\\[5pt]
&=&\frac{1}{(q^3;q^3)_\infty}\sum_{n=0}^\infty q^{n(n+1)/2}.
\end{eqnarray*}
Note that there are no triangular numbers that is congruent to $2$
modulo $3$. It follows that
\[
R_3(0,3,3n+2)= R_3(1,3,3n+2).
\]
Since
$R_3(1,3,3n+2)=R_3(2,3,3n+2)$, the proof is complete. \qed

To conclude this section, we remark that the rank $r_3$ can be
used to give combinatorial explanations of
many Ramanujan-type
congruences for $\overline{pp}(n)$ which
plays an analogous role to the rank introduced by
Bringmann and Lovejoy \cite{Lovejoy08}
for congruences also for overpartition pairs.
To be specific, we
have the following theorem. The proof is similar to the
proof of Bringmann and Lovejoy.  But the rank $r_3$ seems to
simpler.

\begin{thm}
 Let $l$ be an
odd prime, and let $t$ be an odd number
 which is a power of $l$ or
relatively prime to $l$. Then for any positive integer $j$, there
are infinitely many non-nested arithmetic
 progressions $An + B$ such
that
\begin{equation}\label{eqcongrank3}
R_3(r,t,An+B)\equiv 0\ ({\rm mod}\ l^j)
\end{equation}
for any $0\leq r\leq t-1$.
\end{thm}

\pf  Note that  the generating function
for $R_3(s,t,n)$  can be
decomposed into a linear
combination of certain modular forms
similar to the case for $\overline{NN}(r,t,n)$.
Suppose that $t$ is
an odd integer and $0\leq s<t$.
Let $\zeta_t=e^{\frac{2\pi i}{t}}$,
and let $C(s,t)$ be the constant given by
\[
C(s,t)=\frac{1}{t}\sum_{k=0}^{t-1}\zeta_t^{-ks}.
 \] Then we have
\[
C(s,t)+\sum_{n=1}^\infty
R_3(s,t,n)q^n=\frac{1}{t}\sum_{k=0}^{t-1}\zeta_t^{-ks}R_3(\zeta_t^k;q),
\]
where
\[
R_3(z;q)=\frac{(-q;q)_\infty^2}{(qz;q)_\infty(q/z;q)_\infty}.
\]
Observe that $R_3(\zeta_t^k;q)$ differs
 from $R(\zeta_t^k;q)$ (see
Bringmann and Lovejoy \cite[proposition 2.4]{Lovejoy08})
only by a factor
$\frac{4}{(1+\zeta_t^k)(1+\zeta_t^{-k})}$.
Hence the argument of Bringmann and Lovejoy for
\eqref{eqrankcong} can be
carried over to deduce relation \eqref{eqcongrank3}.
This completes the proof.  \qed

\section{Infinite families of congruences modulo $3$ and $5$}

In this section, we  obtain a formula
 for $\overline{pp}(3n)$
modulo $3$ based on
the number of representations of $n$ as a sum of two
squares. We further derive a
formula for  $\overline{pp}(5n)$ modulo $5$ in connection with
the number of representations of $n$ in the
form $x^2+5y^2$.
As consequences, we  give infinite families
of congruences modulo $3$
and $5$.

\begin{thm}\label{thmcong3n}
If the prime factorization of $n$ is given by
\begin{equation}\label{eqfacn1}
n=2^a\prod_{i=1}^rp_i^{v_i}\prod_{j=1}^sq_j^{w_j},
\end{equation}
where $p_i\equiv 1\ ({\rm mod\ }4)$ and $q_j\equiv 3\ ({\rm mod\
}4)$. Then
\begin{equation}\label{eqpp3n}
\overline{pp}(3n)\equiv
(-1)^n\prod_{i=1}^r(1+v_i)\prod_{j=1}^s\frac{1+(-1)^{w_j}}{2}\ ({\rm
mod}\ 3).
\end{equation}
\end{thm}

\pf First, it is easy to see that
\[
\varphi(-q)^3\equiv \varphi(-q^3)\ ({\rm mod\ }3)
\]
and
\[
\varphi(-q)=\varphi(-q^9)+qB(q^3),
\]
where $B(q)$ is a infinite series in $q$ with integer coefficients .
Hence,
\begin{eqnarray*}
\sum_{n=0}^\infty
\overline{pp}(n)q^n&=&\frac{\varphi(-q)}{\varphi(-q)^3}
\equiv\frac{\varphi(-q)}{\varphi(-q^3)}\ ({\rm mod}\ 3)\\[5pt]
&=&\frac{\varphi(-q^9)+qB(q^3)}{\varphi(-q^3)}.
\end{eqnarray*}
Extracting the terms $q^{3n}$ for $n\geq 0$, and replacing $q^3$ by
$q$, we find that
\begin{eqnarray}\label{eqpp3na}
\sum_{n=0}^\infty \overline{pp}(3n)q^n&\equiv&
\frac{\varphi(-q^3)}{\varphi(-q)}
\equiv \varphi(-q)^2\ ({\rm mod\ }3).
\end{eqnarray}
Let $r_2(n)$ denote the number of representations of $n$
as a sum of two squares. So we have
\begin{equation}\label{eqpp3nb}
\varphi(-q)^2 = \sum_{n=0}^\infty (-1)^nr_2(n)q^n.
\end{equation}
From (\ref{eqpp3na}) and (\ref{eqpp3nb}) it follows that
\begin{equation} \label{eqpp3nc}
\overline{pp}(3n)\equiv (-1)^nr_2(n)\ ({\rm mod\ }3).
\end{equation}
Given the prime
factorization of $n$ in the form of
  \eqref{eqfacn1}, it is well
known that
\begin{equation}\label{eqpp3nd}
 r_2(n)=4\prod_{i=1}^r(1+v_i)\prod_{j=1}^s
 \frac{1+(-1)^{w_j}}{2},
\end{equation}
see, for example, Berndt \cite{Berndt06} or Grosswald
\cite{Grosswald84}. Combining
 (\ref{eqpp3nc}) and (\ref{eqpp3nd}),
we get the desired formula (\ref{eqpp3n}). \qed

\begin{thm}\label{thmcong3}
Assume that $p$ is prime with $p\equiv 3\ ({\rm mod\ }4)$, and
$s$ is an integer with $1\leq s< p$. Then
for any $\alpha\geq 0$ and $n\geq 0$, we have
\begin{equation}\label{eqp^a}
\overline{pp}(3p^{2\alpha+1}(pn+s))\equiv 0 \ ({\rm mod\ }3).
\end{equation}
In particular, setting $p=3$, we have
 for any $\alpha \geq 1$ and  $n\geq 0$,
\begin{equation}\label{eq3^a}
\overline{pp}(9^\alpha(3n+1))\equiv 0\ ({\rm mod}\ 3)
\end{equation}
and
\begin{equation}\label{eq3^b}
\overline{pp}(9^\alpha(3n+2))\equiv 0\ ({\rm mod}\ 3).
\end{equation}
\end{thm}

\pf Recall that  $r_2(n)=0$
if and only if there exists a prime
congruent to $3$ modulo $4$ that has an odd exponent in the
canonical factorization of $n$. It can be seen that
\[
r_2(p^{2\alpha+1}(pn+s))=0,
\]
since $p$ is not a factor of $pn+s$. By \eqref{eqpp3nc} we obtain
the congruence relation \eqref{eqp^a}. This completes the proof.\qed

\begin{thm}\label{thmcong5n}
Let $R(n,x^2+5y^2)$ denote the number of
 representations of $n$ by the quadratic form $x^2+5y^2$.
Then for any $n\geq 0$, we have
\begin{equation}\label{eqthmpp5n}
\overline{pp}(5n)\equiv (-1)^nR(n,x^2+5y^2)\ ({\rm mod}\ 5).
\end{equation}
\end{thm}

To prove this theorem, we need a simple property of
a Lambert series on the right-hand side of
 \eqref{eqlem}.

\begin{lem}\label{lemaar}
Let $1\leq r\leq 4$ and let the numbers $a_r(n)$ be given by
\begin{equation}\label{eqlem}
\sum_{n=1}^\infty a_r(n)q^n=\sum_{n=0}^\infty
\frac{q^{5n+r}}{1-(-q)^{5n+r}}.
\end{equation}
 Then we have
 \[
 \sum_{n=1}^\infty a_r(5n)q^{n}=\sum_{n=1}^\infty a_r(n)q^n.
\]
\end{lem}

\pf Let $r_8(n)$ denote  the number of representations of $n$ as a sum
of eight  squares, namely,
\begin{equation}\label{eqpp5a}
\varphi(q)^8=1+\sum_{n=1}^\infty r_8(n)q^n.
\end{equation}
Using the generating function (\ref{eqgenpp2}) for
$\overline{pp}(n)$, we have
\begin{equation}\label{eqpp5b}
\varphi(q)^{10}\sum_{n=0}^\infty
\overline{pp}(n)(-q)^n=\varphi(q)^8.
\end{equation}
In view of relation (\ref{eqpp5a}) and the fact that
\[
\varphi(q^5)\equiv \varphi(q)^5\ ({\rm mod\ }5),
\]
from \eqref{eqpp5b}  we see that
\[
\varphi(q^5)^2\sum_{n=0}^\infty
\overline{pp}(n)(-q)^n\equiv
1+\sum_{n=1}^\infty r_8(n)q^n
\ ({\rm mod\ }5).
\]
Choosing  the terms for which the power of
$q$ is a multiple of $5$, we find that
\begin{equation}\label{eqpp5c}
\varphi(q^5)^2\sum_{n=0}^\infty \overline{pp}(5n)(-q)^{5n}\equiv
1+\sum_{n=1}^\infty r_8(5n)q^{5n} \ ({\rm mod\ }5).
\end{equation}
Replacing $q^5$ by $q$ in \eqref{eqpp5c} gives
\begin{equation}\label{eqpp5d}
\varphi(q)^2\sum_{n=0}^\infty \overline{pp}(5n)(-q)^{n}\equiv
1+\sum_{n=1}^\infty r_8(5n)q^{n} \ ({\rm mod\ }5).
\end{equation}
We wish to establish the
following congruence
\begin{equation}\label{eqr5n}
1+\sum_{n=1}^\infty r_8(5n)q^{n}\equiv 1+\sum_{n=1}^\infty r_8(n)q^n
\ ({\rm mod\ }5).
\end{equation}
To this end, we recall that
\begin{equation}\label{eqsquare8}
\varphi(q)^8=1+16\sum_{n=1}^\infty \frac{n^3q^n}{1-(-q)^n},
\end{equation}
see Berndt \cite[Theorem 3.5.3]{Berndt06}.
It follows that
\begin{eqnarray}
\varphi(q)^8&\equiv&1+\sum_{n=0}^\infty\frac{q^{5n+1}}{1-(-q)^{5n+1}}+
\sum_{n=0}^\infty\frac{3q^{5n+2}}{1-(-q)^{5n+2}}
\nonumber \\[5pt]
& & \;\;+\sum_{n=0}^\infty\frac{2q^{5n+3}}{1-(-q)^{5n+3}}+
\sum_{n=0}^\infty\frac{4q^{5n+4}}{1-(-q)^{5n+4}}\
 ({\rm mod\ }5). \label{eqsquare8a}
\end{eqnarray}
For $r=1,2,3,4$, let $a_r(n)$ be defined by (\ref{eqlem}). By
(\ref{eqpp5a}), the above relation (\ref{eqsquare8a}) can be
rewritten as
\begin{equation}\label{eqsquare8b}
\sum_{n=1}^\infty r_8(n)q^n\equiv \sum_{n=1}^\infty
a_1(n)q^n+3\sum_{n=1}^\infty a_2(n)q^n+2\sum_{n=1}^\infty
a_3(n)q^n+4\sum_{n=1}^\infty a_4(n)q^n\ ({\rm mod\ }5).
\end{equation}
By Lemma \ref{lemaar}, we find that
\begin{eqnarray}
\sum_{n=1}^\infty r_8(5n)q^{n}&\equiv& \sum_{n=1}^\infty
a_1(5n)q^n+3\sum_{n=1}^\infty a_2(5n)q^n+2\sum_{n=1}^\infty
a_3(5n)q^n+4\sum_{n=1}^\infty a_4(5n)q^n\ ({\rm mod\ }5)
\nonumber\\[5pt]
&=&\sum_{n=1}^\infty a_1(n)q^n+3\sum_{n=1}^\infty
a_2(n)q^n+2\sum_{n=1}^\infty a_3(n)q^n+4\sum_{n=1}^\infty
a_4(n)q^n.\label{eqsquare8c}
\end{eqnarray}
By \eqref{eqsquare8b} and \eqref{eqsquare8c} we obtain
\eqref{eqr5n}. From \eqref{eqpp5d} and \eqref{eqr5n}  we know that
\[
\varphi(q)^2\sum_{n=0}^\infty \overline{pp}(5n)(-q)^{n} \equiv
\varphi(q)^8 \ ({\rm mod\ }5).
\]
Thus,
\begin{equation}\label{eqpp5e}
\sum_{n=0}^\infty \overline{pp}(5n)(-q)^{n}
\equiv \varphi(q)^6 \equiv \varphi(q)\varphi(q^5)
 \ ({\rm mod\ }5).
\end{equation}
By the definition of $R(n,x^2+5y^2)$, we see that
\begin{equation}\label{eqpp5f}
\sum_{n=0}^\infty R(n,x^2+5y^2)q^n=\varphi(q)\varphi(q^5).
\end{equation}
As a consequence of \eqref{eqpp5e} and \eqref{eqpp5f},
 we deduce that
\[
\sum_{n=0}^\infty \overline{pp}(5n)(-q)^{n}\equiv
\sum_{n=0}^\infty R(n,x^2+5y^2)q^n \ ({\rm mod\ }5).
\]
Thus the proof is complete by equating  coefficients. \qed

The formula for $R(n,x^2+5y^2)$ due to
 Berkovich and Yesilyurt
\cite{berkovich09} leads to the following
formula for
$\overline{pp}(5n)$ modulo $5$.

\begin{thm}If the prime factorization of $n$ is given by
\begin{equation}\label{eqfacn2}
n=2^a5^b\prod_{i=1}^rp_i^{v_i} \prod_{j=1}^s q_j^{w_j},
\end{equation}
where $p_i\equiv 1,3,7$, or $9\ ({\rm mod\ }20)$  and $q_j\equiv
11,13,17$, or $19\ ({\rm mod\ }20)$. Then we have
\begin{equation}\label{eqcong5n}
\overline{pp}(5n)\equiv (-1)^n
\left(1+(-1)^{a+t}\right)\prod_{i=1}^r(1+v_i)\prod_{j=1}^s
\frac{1+(-1)^{w_j}}{2}\ ({\rm mod}\ 5),
\end{equation}
where $t$ is the number of prime factors of $n$, counting
multiplicity, that are congruent to $3$ or $7$ modulo $20$.
\end{thm}

\pf Given the prime factorization of $n$ in the form of
\eqref{eqfacn2}, it is known that
\begin{equation}\label{eqR}
R(n,x^2+5y^2)=\left(1+(-1)^{a+t}\right)\prod_{i=1}^r(1+v_i)\prod_{j=1}^s
\frac{1+(-1)^{w_j}}{2},
\end{equation}
 see, for example, Berkovich and
Yesilyurt \cite[Corollary 3.3]{berkovich09}. Combining
\eqref{eqthmpp5n} and (\ref{eqR}), we get (\ref{eqcong5n}).
This completes the proof. \qed

Based on the above theorem, we establish two infinite families of
congruences modulo 5.

\begin{thm}\label{thmcong5^a}
For any $\alpha \geq 1$ and  $n\geq 0$, we have
\begin{equation}\label{eqcong5^a1}
\overline{pp}(5^\alpha(5n+2))\equiv 0\ ({\rm mod}\ 5)
\end{equation}
and
\begin{equation}\label{eqcong5^a2}
\overline{pp}(5^\alpha(5n+3))\equiv 0\ ({\rm mod}\ 5).
\end{equation}
\end{thm}

\pf Considering the possible residues of $x^2+5y^2$
modulo $5$, we
find that
\[
R(5n+2,x^2+5y^2)=R(5n+3,x^2+5y^2)=0.
\]
In light of \eqref{eqthmpp5n}, we deduce that
\begin{equation}\label{eqR5n+2}
\overline{pp}(25n+10) \equiv (-1)^{5n+2}R(5n+2,x^2+5y^2)
\equiv 0 \ ({\rm
mod}\ 5)
\end{equation}
and
\begin{equation}\label{eqR5n+3}
\overline{pp}(25n+15) \equiv (-1)^{5n+3}R(5n+3,x^2+5y^2)
\equiv 0 \ ({\rm
mod}\ 5).
\end{equation}
Observe that   formula
\eqref{eqcong5n} for $\overline{pp}(5n)$
 modulo $5$ is independent of the exponent
 of $5$ in the
factorization of $n$. This means  that
 for $\alpha\geq 1$,
\begin{equation}\label{eqcong5^a3}
\overline{pp}(5n)\equiv
\overline{pp}(5^\alpha n) \ ({\rm mod}\ 5).
\end{equation}
 Combining \eqref{eqR5n+2}, \eqref{eqR5n+3}  and
 \eqref{eqcong5^a3}, we obtain the desired congruence relations (\ref{eqcong5^a1}) and
 (\ref{eqcong5^a2}).  This completes
 the proof.
 \qed

\section{Further congruences for overpartition pairs }

In this section, we shall establish some congruences for
$\overline{pp}(n)$ modulo 5 and 9 which are similar to the
congruences for the number of broken $2$-diamonds partitions derived
by Paule and Radu \cite{Paule09}. Let us begin with the congruences
modulo 9 which are derived from congruences modulo 3.

\begin{thm}\label{thmstrange1}
For any prime with $p\equiv 1\ ({\rm mod\ 12})$, we have
\begin{equation}
\overline{pp}((3n+2)p)\equiv
\frac{\overline{pp}(2p)}{3}\overline{pp}(3n+2)\ ({\rm mod\ 9}),
\end{equation}
for all positive integers $n$ such that $3n+2\not\equiv 0\ ({\rm
mod\ }p)$.
\end{thm}

To prove the above theorem, we need the following lemma
which is a
special case of  Newman \cite[Theorem 3]{Newman59}.

\begin{lem}\label{lemnewman1}
For each prime $p$ with $p\equiv 1\ ({\rm mod\ }12)$
 and  for all
positive integers $n$,
\begin{equation}\label{eqnewman1}
b\left(np+\frac{2p-2}{3}\right)+
p^4b\left(\frac{n}{p}-2\frac{p-1}{3p}\right)=
b\left(\frac{2p-2}{3}\right)b(n),
\end{equation}
where $b(n)$ is defined by
\[
\sum_{n=0}^\infty b(n)q^n=(q;q)_\infty^4(q^2;q^2)_\infty^6.
\]
\end{lem}

Since the equality is derived by equating coefficients
of series in $q$, it is safe to assume that
 $b(t)=0$  if $t$ is not a
nonnegative integer. We are now ready to give a proof of Theorem
\ref{thmstrange1}.

\noindent {\it Proof of Theorem \ref{thmstrange1}.} By the
generating function of $\overline{pp}(3n+2)$ as given in
\eqref{eq3n+2}, we see that
\[
\sum_{n=0}^\infty \frac{\overline{pp}(3n+2)}{3}q^n\equiv
(q;q)_\infty^4(q^2;q^2)_\infty^6\ ({\rm mod\ }3).
\]
From the definition of $b(n)$, we deduce that for $n\geq 0$,
\begin{equation}\label{eqstrange1a}
\frac{\overline{pp}(3n+2)}{3}\equiv b(n)\ ({\rm mod\ }3).
\end{equation}
On the other hand, for those prime $p$ with
$p\equiv 1\ ({\rm mod\ }12)$ and those $n$
such that $3n+2$ is not a
multiple of $p$, it follows that
$b\left(\frac{n}{p}-2\frac{p-1}{3p}\right)=0$. Thus, by
Lemma \ref{lemnewman1} we obtain
\begin{equation}\label{eqstrange1b}
b\left(np+\frac{2p-2}{3}\right)= b\left(\frac{2p-2}{3}\right)b(n).
\end{equation}
Substituting \eqref{eqstrange1a} into \eqref{eqstrange1b}, we get
\[
\frac{1}{3}\overline{pp}(3np+2p)
\equiv \frac{1}{9}\overline{pp}(2p)
\overline{pp}(3n+2)\ ({\rm mod\
}3),
\]
as required.\qed

Next, we use Lemma \ref{lemnewman1} to
obtain the following congruences. Such type of
congruences are called strange congruences
by Paule and Radu \cite{Paule09}.

\begin{thm}\label{thmstrange2} For any
$k\geq 0$, we have
\begin{equation}\label{eqstrange2a}
\overline{pp}(2\cdot 13^k)\equiv 3(k+1)\ ({\rm mod\ }9).
\end{equation}
\end{thm}

\pf Let $p$ be a prime with $p\equiv 1\ ({\rm mod\ }12)$. Then
 setting $n=2(p^{k+1}-1)/3$ in
\eqref{eqnewman1} and using \eqref{eqstrange1a}, we get
\[
\frac{1}{3}\overline{pp}(2p^{k+2})+\frac{1}{3}\overline{pp}(2p^k)\equiv
\frac{1}{9}\overline{pp}(2p)\overline{pp}(2p^{k+1})\ ({\rm mod\ }3).
\]
When  $p=13$,
since $\overline{pp}(26)\equiv 6\ ({\rm mod\ }9)$, we
deduce that
\begin{equation}\label{eqcong13^a}
\overline{pp}(2\cdot 13^{k+2})+\overline{pp}(2\cdot 13^k)\equiv
2\overline{pp}(2\cdot 13^{k+1})\ ({\rm mod\ }9).
\end{equation}
Given the initial conditions
 $\overline{pp}(2)\equiv 3\ ({\rm mod\ }9)$,
 $\overline{pp}(26)\equiv 6\ ({\rm mod\ }9)$, by
 iteration of  \eqref{eqcong13^a}, we
reach the conclusion \eqref{eqstrange2a}. This completes the proof.
\qed

We now turn to the congruences modulo $5$.

\begin{thm}\label{thmstrange3} Let
$p$ be a  prime with $p\equiv 13\ ({\rm mod\ 20})$ or
 $p\equiv 17\ ({\rm mod\ 20})$. Then the following congruence
 holds any positive integer  $n$  that is not divisible by $p$,
 \begin{equation}
 \overline{pp}(5pn)\equiv 0\ ({\rm mod\ }5).
 \end{equation}
\end{thm}

The following lemma is a special case of Newman \cite[Theorem
3]{Newman59}, which will be needed in the proof of
 Theorem \ref{thmstrange3}.

\begin{lem}\label{lemnewman2}
For each prime $p$ with $p\equiv 1\ ({\rm mod\ }4)$
 and  for all
positive integers $n$,
\begin{equation}\label{eqnewman2}
c(np)+p^2c(n/p)=(p^2+1)c(n),
\end{equation}
where $c(n)$ is defined by
\[
\sum_{n=0}^\infty
c(n)q^n=\varphi(-q)^6=\frac{(q;q)_\infty^{12}}{(q^2;q^2)_\infty^6}.
\]
\end{lem}

Like the case for
  Lemma \ref{lemnewman1}, we may assume that
 $c(t)=0$  if $t$ is not a
nonnegative integer.

\noindent {\it Proof of Theorem \ref{thmstrange3}.} Recall the
following relation as given in  \eqref{eqpp5e},
\[
\sum_{n=0}^\infty \overline{pp}(5n)q^n \equiv \varphi(-q)^6 \ ({\rm
mod\ }5).
\]
From the definition of $c(n)$, we see that for $n\geq 0$,
\begin{equation}\label{eqstrange3a}
\overline{pp}(5n)\equiv c(n)\ ({\rm mod\ }5).
\end{equation}
On the other hand, for any prime $p$
with $p\equiv 13\ ({\rm mod\
20})$ or
 $p\equiv 17\ ({\rm mod\ 20})$, and for any
 $n$  that is not a multiple of $p$, we have
$c(n/p)=0$. Thus, by Lemma \ref{lemnewman2} we obtain that
\begin{equation}\label{eqstrange3b}
c(np)= (p^2+1)c(n).
\end{equation}
Combining \eqref{eqstrange3a} and \eqref{eqstrange3b}, we get
\begin{equation}\label{eqstrange3c}
\overline{pp}(5np)\equiv (p^2+1)\overline{pp}(5n)\ ({\rm mod\ }5).
\end{equation}
Since $p^2+1$ is a multiple of $5$ if  $p$ is a prime
of the form $20k+13$ or
$20k+17$, the above congruence reduces to
\[
\overline{pp}(5np)\equiv 0\ ({\rm mod\ }5).
\]
 This
completes the proof.\qed

To conclude this section, we use Lemma \ref{lemnewman2} to
derive the following congruences.

\begin{thm}\label{thmstrange4}  Let $p$ be a
prime with $p\equiv 1\ ({\rm mod\ }20)$ or
 $p\equiv 9\ ({\rm mod\ }20)$. Then the following
 congruence holds for any positive integer $k$,
 \begin{equation}
 \overline{pp}(5p^k)\equiv 3(k+1)\ ({\rm mod\ }5).
 \end{equation}
\end{thm}

\pf
 Applying \eqref{eqnewman2} with $n=p^{k+1}$
 and using \eqref{eqstrange3a},
we get
\begin{equation}\label{eqstrange4a}
\overline{pp}(5p^{k+2})+p^2\overline{pp}(5p^k)\equiv
(p^2+1)\overline{pp}(5p^{k+1})\ ({\rm mod\ }5).
\end{equation}
It is easily seen that $p^2\equiv 1 \ ({\rm mod\ }5)$. From
\eqref{eqstrange4a} we see that
\begin{equation}\label{eqstrange4b}
\overline{pp}(5p^{k+2})+\overline{pp}(5p^k)\equiv
2\overline{pp}(5p^{k+1})\ ({\rm mod\ }5).
\end{equation}
Setting $n=1$ in \eqref{eqstrange3c}, we get
\[
\overline{pp}(5p)\equiv (p^2+1)\overline{pp}(5)\equiv
2\overline{pp}(5)\ ({\rm mod\ }5).
\]
Since $\overline{pp}(5)\equiv 3\ ({\rm mod\ }5)$,
so
$\overline{pp}(5p)\equiv 1\ ({\rm mod\ }5)$. By iteration
 of \eqref{eqstrange4b}, we  arrive at
 the desired congruence, and hence the proof is complete.
\qed

\vspace{.2cm} \noindent{\bf Acknowledgments.} This work was
supported by the 973 Project, the PCSIRT Project of the Ministry of
Education, and the National Science Foundation of China.

\end{document}